\documentclass{article}

\newfont{\sheaf}{eusm10 scaled\magstep1} 
 
\newcommand{\hol}{\ensuremath{\mathcal{O}}} 
\def\Bbb{\bf}
\def\PP{{\Bbb P}}
\def\R{{\Bbb R}}
\def\C{{\Bbb C}}
\def\Q{{\Bbb Q}}
\def\Z{{\Bbb Z}}

\newtheorem{teo}{Theorem}[section] 
\newtheorem{df}[teo]{Definition} 
\newtheorem{lem}[teo]{Lemma} 
\newtheorem{cor}[teo]{Corollary} 
\newtheorem{oss}[teo]{Remark} 
\newtheorem{prop}[teo]{Proposition} 

\newenvironment{proof}{\medskip {\bf Proof.}}{\hfill \rule{.5em}{1em}
\\}

\def\eea{\end{eqnarray*}}
\def\bea{\begin{eqnarray*}}

\begin{document}

\title{Symplectic structures of algebraic surfaces and deformation.}

\author{Fabrizio Catanese\\
 Universit\"at Bayreuth \thanks{
The research of the  author was performed in the realm  of the 
 SCHWERPUNKT "Globale Methode in der komplexen Geometrie",
and of the EAGER EEC Project. 
}
 \\
This article is dedicated to the memory of Boris  Moisezon }

\date{May 18, 2002}
\maketitle

\begin{abstract} 
 We show that a surface of general type has a canonical 
symplectic structure (up to symplectomorphism)
 which is  invariant for smooth deformation. 
 Our main theorem is that the symplectomorphism type is also invariant
 for deformations which also allow certain normal singularities, called
 Single Smoothing Singularities ( and abbreviated as SSS),
or yielding $\Q$-Gorenstein smoothings of quotient singularities. 
 
Using the counterexamples of M.Manetti to the 
 DEF = DIFF question whether deformation type and 
diffeomorphism type  coincide for algebraic surfaces,
we show that  these yield surfaces of general type which are 
not deformation equivalent but are symplectomorphic.
In particular, they are diffeomorphic through a diffeomorphism
carrying the canonical class of one to the canonical class of
the other surface.

ii)  Another interesting corollary is the existence of cuspidal 
algebraic plane curves which are symplectically isotopic, but not
equisingular deformation equivalent. 

  \end{abstract}

\vfill
\pagebreak
\section{Introduction}

The present note is an addendum to the paper 
"Moduli spaces and real structures" (\cite{cat4}), where simple examples were
 shown   of surfaces which are diffeomorphic, but not
deformation equivalent (in other words,  belonging to different
 connected components of the moduli space).
For those examples, one would take complex conjugate surfaces,
whence the condition of being orientedly diffeomorphic was tautologically
fulfilled. However, the diffeomorphism would send the canonical class
to its opposite, and one could also rephrase the heart of the proof
as saying that these surfaces admit no self-homeomorphism reversing
the canonical class.

Marco Manetti ( \cite{man4} had earlier constructed examples 
of surfaces of
general type which are not deformation equivalent, and which admit
a common degeneration to a normal surface with singularities, yielding
each a $\Q$-Gorenstein smoothing of these singularities.
Using a result of Bonahon on the group of diffeomorphisms 
of lens spaces he was able to show that the surfaces are diffeomorphic.

Indeed, with a direct proof a more general result holds true:
 we admit deformations with singular fibres which are normal and 
 in the following class SSS

\begin{df}
A Single Smoothing Surface Singularity $(X_0, x_0)$ is a normal 
surface singularity such that the smoothing locus
$$\Sigma  := \{ t \in Def(X_0, x_0) | X_t \ is \ smooth \} $$
is irreducible (cf. \cite{l-w} and references therein for
examples of S.S.S. singularities, such as complete intersections,
some cusp and triangle singularities..).
 
\end{df}

\begin{teo}
Let $  \mathcal {X}\subset \PP^N \times \Delta$ and 
$  \mathcal {X}' \subset \PP^N \times \Delta'$ be two flat families
of normal surfaces over the disc of radius 2 in $\C$ such that

1) the central fibres are projectively equivalent in  $\PP^N$ , and 
$X_0 = X_0'$ is a surface with single smoothing singularities 

2) the other fibres $X_t$, $X_t'$, for $ t, t' \neq 0$ are smooth.

Set $X := X_1$, $X' := X_1'$: then

a) $X$ and $X'$ are diffeomorphic

b) if  $FS $ denotes the symplectic form inherited from the Fubini-Study
K\"ahler metric on $\PP^N$, then the symplectic manifolds
$(X,FS)$ and $(X',FS)$ are symplectomorphic.

The same conclusion holds if  hypothesis 1) is replaced by 

1') the two flat families yield
  $\Q$-Gorenstein smoothings of the singularities of $X_0 = X'_0$.
\end{teo}

A corollary of the proof of the above theorem is the following

\begin{teo} A minimal surface of general type $S$ has a canonical symplectic
 structure, unique up to symplectomorphism, such that the class of
the symplectic form is the class of the canonical sheaf
$ \omega_S =  \Omega^2_S = \hol_S (K_S).$

\end{teo}

\begin{oss}
Theorem 1.2 holds more generally under the assumption that 
the two smoothings $  \mathcal {X},   \mathcal {X}'$ of $X_0$ lie, 
for each singular point $ x_0$ of $X_0$, in the same irreducible 
component of $Def (X_0, x_0) $.
\end{oss}

As already mentioned, for the main application 
we need to apply Theorem 1.2 under hypothesis 1')
above. 

 To this purpose we recall some known facts on the
 class of singularities given by the (cyclic)
quotient singularities admitting a $\Q$-Gorenstein smoothing
(cf. \cite{man4}, Section 1, pages 34-35, or the original sources
\cite{k-sb},\cite{man0},\cite{l-w}).

The simplest way to describe these singularities  
  $${\bf Cyclic \ quotient \ singularity} \ \frac{1}{d n^2} (1,dna-1) = A_{dn-1} / \mu_n$$
is to view them on the one side as quotients of $\C^2$ by a cyclic 
group of order $d n^2$ acting with the indicated characters $(1,dna-1)$,         
or as quotients  of the rational double point $ A_{dn-1}  $
of equation $ uv - z^{dn} = 0$ by the action of the group $\mu_n$ of
n-roots of unity acting in the following way:
 
$$\xi \in \mu_n  \ {\bf acts \  by \ :} \ 
 (u,v,z) \rightarrow  ( \xi u , \xi^{-1} v,\xi^{a} z).$$

This quotient action gives rise to a quotient family 
$ \mathcal{X} \rightarrow \C^d$, where 

$ \mathcal{X}= 
 \mathcal{Y}/ \mu_n$ , $ \mathcal{Y}$ is the hypersurface in
$\C^3 \times \C^d$ of equation
$ uv - z^{dn} = \Sigma_{k=0}^{d-1} t_k z ^{kn}$ and the action of 
$\mu_n$ is extended trivially on the factor $\C^d$.

The heart of the construction is that $\mathcal{Y}$, being
 a hypersurface,
is Gorenstein (this means that the canonical sheaf
 $\omega_{\mathcal{Y}}$ is 
invertible), whence such a quotient $ \mathcal{X}= 
 \mathcal{Y}/ \mu_n$, by an action which is unramified in codimension $1$,
is (by definition) $\Q$-Gorenstein.

These smoothings were considered by Kollar and Shepherd Barron 
in \cite{k-sb} (and independently, Manetti,
 \cite{man0}) who showed their relevance in the theory
of compactifications of moduli spaces of surfaces.

Riemennschneider ( \cite{riem}) earlier showed that for these cyclic quotient 
singularities the basis of the
semiuniversal deformation  can consist of two smooth components
 crossing transversally, each one yielding a smoothing, but
only one admitting a simultaneous resolution, and only the other
yielding  smoothings with $\Q$-Gorenstein total space.

Recently Marco Manetti ( \cite{man4}) was able to produce examples
of surfaces of general type which are not deformation equivalent, 
but diffeomorphic. 

His examples are based on a complicated construction  of Abelian 
coverings  of rational surfaces with group $( \Z /2)^m$, leading to
families $ \mathcal{X}$, $ \mathcal{X}'$ as in Theorem 1.2, since
  the two families induce smoothings of cyclic quotient singularities
which are $\Q$-Gorenstein. Whence follows right away

\begin{teo}
Manetti's surfaces (Section 6 in \cite{man4}) 
provide examples of surfaces of general type
which are not deformation equivalent, but, endowed with their 
canonical symplectic structures, are symplectomorphic.
\end{teo} 

Manetti's surfaces are, like our examples (\cite{cat4}) and the ones 
by Kharlamov-Kulikov (\cite{k-k}), not simply connected. 

If one would insist on finding counterexamples to the Friedman-Morgan
conjecture which are simply connected, there are several  
 candidates. Some very natural ones are illustrated  in the next 
section, where we shall show several families
of surfaces which are  pairwise homeomorphic by a homeomorphism
 carrying the canonical 
class to the canonical class, but are not deformation equivalent.

Finding pairs of families yielding diffeomorphic surfaces would give 
the simply connected counterexamples,
finding a pair of families yielding non diffeomorphic surfaces would 
be even more interesting.

\section{A digression and a problem.}

This section proposes two examples 
(the second depends on 3 integer parameters (a,b,c)) of  families of surfaces,
 which are homeomorphic by a homeomorphism carrying the canonical class
to the canonical class, in view of the following  proposition 
(well known to experts) 

\begin{prop}
Let $S$, $S'$ be simply conected minimal surfaces of general type
such that $p_g(S) = p_g(S') \geq 1$ , $K^2_S = K^2_{S'}$, and moreover such
 that the divisibility indices of $K_S $ and $K_{S'}$ are the same.

Then there exists a homeomorphism $F$ between $S $ and $S'$, 
unique up to isotopy,
carrying $K_{S'}$ to $K_S $.
\end{prop}

\begin{proof}
By Freedman's theorem (\cite{free}, cf. especially 
\cite{f-q}, page 162) for each isometry
$h: H_2(S, \Z) \rightarrow  H_2(S', \Z)$ there exists a homeomorphism 
$F$ between $S $ and $S'$, 
unique up to isotopy, such that $F_{*} = h$. In fact, $S$ and $S'$ are
smooth $4$-manifolds, whence
 the Kirby-Siebenmann invariant vanishes.

Our hypotheses that $p_g(S) = p_g(S')$ , $K^2_S = K^2_{S'}$
and that $K_S, K_{S'}$ have the same divisibility
imply that the two lattices $ H_2(S, \Z)$, $ H_2(S', \Z)$ have the
same rank, signature and parity, whence they are isometric
(since $S, S'$ are algebraic surfaces, 
cf. e.g. \cite{cat1}). Finally, by Wall's theorem (\cite{wall}) (cf. also
\cite{man2}, page 93) such isometry $h$ exists since the vectors
corresponding to the respective canonical classes have the same 
divisibility and by Wu's theorem they are characteristic:
 in fact Wall's condition
$ b_2 - |\sigma| \geq 4$ 
($\sigma$ being the signature of the intersection form) 
is equivalent to $p_g \geq 1$. 
\end{proof}

For the families in the second example we show that they yield different
deformation types. In both examples we are unable yet to decide
whether the  surfaces belonging to different families 
are diffeomorphic to each other.

As in \cite [Sections 2,3,4]{cat1} we consider smooth bidouble covers 
$S$ of  ${\Bbb P^1\times 
\Bbb P^1}$: these are smooth finite Galois covers of  ${\Bbb P^1\times 
\Bbb P^1}$ having Galois group
 $({\Bbb Z}/2)^2.$ 
Bidouble covers are divided into those of simple type , and those
 not of simple type . 

Those of {\em simple type} (and type $(2a,2b),(2c,2d)$) are defined by 2
equations 
\begin{eqnarray}  z^2 &=&
 f(x,y)\\
 w^2 &=&
 g(x,y) \nonumber ,\end{eqnarray}

where f and g are bihomogeneous polynomials ,  belonging to
 respective vector spaces of sections of line bundles: 
 $H^0({\Bbb P}^1 \times {\Bbb P}^1, {\cal O}(2a,2b)) $ and 
  
$H^0({\Bbb P}^1\times {\Bbb P}^1, {\cal O}(2c,2d)).$
\bigskip

In general, (cf.\cite{cat1}) bidouble covers are embedded in the total space 
of the direct sum of the inverses of 3 line bundles $L_i$ , and 
defined there by equations
( where (i,j,k) is a permutation of {1,2,3} )
\begin{eqnarray}  z_i^2 &=&
 f_j(x,y) f_k(x,y)\\
 z_k f_k(x,y)  &=&
 z_i   z_j \nonumber .\end{eqnarray}

If there is an index $k$ with $f_k(x,y) = 1$ , 
then we have a cover of simple type.
 Therefore we define {\em non simple type} when all the 3 branch divisors
$D_j = \{ f_j(x,y) = 0 \} $ are non trivial . 

Notice moreover that the smoothness of $S$ is ensured 
by the condition that the 3 branch
 divisors be smooth and intersect transversally .

If $(n_i,m_i)$ are the bidegrees of $D_i$ , we shall say that $S$
 is of type

 $(n_1,m_1),(n_2,m_2),(n_3,m_3))$

and we observe that  permuting the indices or exchanging the
$n_i$'s with $m_i$'s does not change the type. 

Moreover,  the $n_i$ 's  are 
either all even or all odd , and likewise  the $m_i$ 's .

In fact  
 $${\cal O}( L_i) = {\cal O}(\frac{1} {2} (n_j + n_k),
\frac{1} {2}(m_j + m_k)).$$

We recall from  \cite [Sections 2,3,4]{cat1} that our surface $S$ has
the following invariants : 
\begin{eqnarray}
 setting \ \  n = \Sigma_{j=1}^3 n_j, 
 m =\Sigma_{j=1}^3 m_j,   \label{symm}\\
\chi (S)= \frac{1}{4} ((n-4)(m-4) + \Sigma_{j=1}^3 n_j m_j)  
 \label{symm}\\
  K^2(S) = 2 (n-4) (m-4)  \label{symm} .\end{eqnarray}

Moreover,(cf.  \cite [Proposition 2.7]{cat1} assume we have a  bidouble cover
 where every branch divisor
$D_j$ is either trivial or  has $n_j \geq 1, m_j \geq 1$.
 Then our surface $S$ is simply connected  unless we have
an "even"   non simple cover 
( this means that there is no trivial $D_j$ , and  all the
 $n_j , m_j $ are divisible by  $2$ : in this case 
% For an "even" non simple cover , 
the fundamental group is $({\Bbb Z}/2)$). 

In order to calculate the divisibility index of the canonical class of such
covers, we will use lemma 4 of \cite{cat3}   for the second class
 of examples.

\begin{df}
Example  1  consists of two  simple covers $S$, $S'$ of respective types
((5,2),(3,2),(1,2)) and ((3,2),(3,2),(3,2)) .
By the above formulae , these two surfaces have the same invariants 
$\chi (S)= 7 , K^2(S) = 20 $. 
\end{df}

\begin{oss}
Clearly, the divisibility of $K$ is at most $2$. However, if
$K$ would be $2$-divisible, then  we would have $ 8 | K^2$,
since if $K \equiv 2 L$, then $ L^2 \equiv L K (mod \ 2) 
\equiv  2 L^2 (mod \ 2) \equiv 0 (mod \ 2) $.
Therefore $K$ is indivisible. 
\end{oss}

\bigskip

The above two families of surfaces  are constructed according to the
same trick used for the next example 2. 
Although we have not yet been able to see whether the two families
 yield diffeomorphic surfaces, nor do we have yet a proof that the two
families are not deformation equivalent, we have decided 
to show this example because of the very low values of 
the numerical invariants.

The families however belong to different irreducible components
 of the moduli space
since the calculations of  \cite [p. 500]{cat1}  show that for the first
 surfaces  the local dimension of the moduli space  is 
at least 39 , whereas for the second ones the local moduli 
pace is smooth of
dimension equal to 38 . 

The above calculations are based on the concept of 
{\em natural deformations
of bidouble covers } (ibidem , def. 2.8 , page 494 )
which will be also used in the forthcoming Theorem 2.6. 

Natural deformations are parametrized by bihomogeneous polynomials 
$f_j(x,y)$ of bidegree
$(n_j,m_j) , \phi_j(x,y)$ of bidegree $(n_j - \frac{1} {2} (n_i + n_k),
m_j - \frac{1} {2} (m_i + m_k))$ , and give equations 
\begin{eqnarray}  z_i^2 &=&   ( f_j(x,y) +\phi_j(x,y) z_j
 ) ( f_k(x,y) + \phi_k(x,y) z_k )\\
  z_i   z_j &=& z_k f_k(x,y) + \phi_k(x,y) z_k^2 
 \nonumber .\end{eqnarray}

In the case of simple covers these specialize to
\begin{eqnarray}  z^2 &=&
 f(x,y) + w \phi(x,y)\\
 w^2 &=&
 g(x,y) + z \psi(x,y) \nonumber ,\end{eqnarray}

and $\phi$ has bidegree (2a-c,2b-d) , whereas $\psi$ has
bidegree (2c-a,2d-b) .

\begin{df}
Example $(a,b,c) $ consists of two  simple covers $S$, $S'$ of respective
 types $((2a, 2b),(2c,2b)$, and $(2a + 2, 2b),(2c - 2,2b)$.
We shall moreover assume, for technical reasons , that $ a \geq 2c + 1$,
 $ a \geq b+2$ ,  $ c \geq b+2$ , and that a,b,c are even and $\geq 3$.
 
By the previous formulae , these two surfaces have the same invariants 
$\chi (S)= 2 (a+c-2) (b-1) + 4 b (a+c) ,  K^2_{S} = 16 (a+c-2) (b-1) $ . 
\end{df}

\begin{oss}
The divisibility index of the canonical divisor $K$ for the above family
of surfaces is easily calculated by lemma 4 of \cite{cat3}, asserting
that the pull back of $H^2(\PP^1 \times \PP^1, \Z)$ is primitively
embedded in $H^2(S, \Z)$. Now, $K_S$ is the pull back of a divisor
of bidegree $(a+c-2,2 b-2)$ whence its divisibility equals simply
$ G.C.D. (a+c-2,2 b-2)$. Therefore the divisibility index is the same
for the several families (vary the integer $k$) of simple covers of 
 types  $(2a + 2k, 2b),(2c - 2k,2b)$.
\end{oss}

\medskip

In the case of example (a,b,c) the natural deformations of S ( which 
yield  all the local deformations , by (ibidem , thm 3.8) ) do not 
preserve the action of the Galois group $({\Bbb Z}/2)^2$ ( this
 would be the case for a cover of type (2a, 2b),(2c,2d) with 
 $ a \geq 2c + 1$ ,  $ d \geq 2b + 1$ , cf.  \cite {cat1,cat2}).

But, since the natural deformations yield equations of type 
\begin{eqnarray}  z^2 &=&
 f(x,y) + w \phi(x,y)\\
 w^2 &=&
 g(x,y)  \nonumber ,\end{eqnarray}
there is preserved the  $({\Bbb Z}/2)$ action sending 
$$ (z,w) \rightarrow (-z,  w).$$
and also the action of  $({\Bbb Z}/2)$ on the quotient 
of $S$ by the above 
involution. 

That is , every small deformation preserves the structure of 
{\em iterated double cover } ( \cite {man3, man4}) . 

We prove now the main result of this section  

\begin{teo}Let  $S$, $S'$ be simple bidouble covers of  ${\Bbb P}^1
\times  {\Bbb P}^1  $ of respective
 types ((2a, 2b),(2c,2b), and (2a + 2k, 2b),(2c - 2k,2b) , where $a,b,c$
are strictly  even integers $\geq 4$, with $ a \geq 2c + 1$,
 $ a \geq b+2$ ,  $ c \geq b+2$ , $ c \geq k+4$ . 
Then  $S$ and  $S'$ are not deformation equivalent.
\end{teo}

\begin{proof}

By \cite [Thm. 3.10] {man4} ,  since conditions C1),C2),
C4), C5) and the first half of condition C3)
are fulfilled , it follows that the connected component of the
 moduli space containing the point corresponding to $S$ consists 
of iterated
double covers of  ${\Bbb P}^1
\times  {\Bbb P}^1  $, 
or of the Segre-Hirzebruch surfaces ${\Bbb F}_{2k}  $, 
which are the deformations of $S$ given by similar formulae
to formula (8).
Rerunning all the arguments in \cite{man4} we notice that
they work verbatim for the more general class of the surfaces 
which are iterated double covers of the Segre-Hirzebruch
 surfaces ${\Bbb F}_{2k}  $.
\end{proof}

We finally observe

\begin{oss}
Consider the two families of surfaces of example $(a,b,c)$ ( def. 2.4),
and let $S$, $S'$ belong to each one of the respective families.

Then there exist two $4$-manifolds with boundary $M_1$, $M_2$ such
that $S$ and $S'$ are obtained by glueing $M_1$ and $M_2$ through 
two  respective 
glueing maps $\phi, \phi' \in Diff (\partial M_1, \partial M_2)$.

These can be described easily as follows: cut $\PP^1$ into two
disks $\Delta_0, \Delta_{\infty}$ and write $\PP^1 \times\PP^1$
as $(\Delta_0 \times \PP^1) \cup  (\Delta_{\infty} \times \PP^1)$.

We let $D'_1$ be the union of $2b$ horizontal lines in
 $\PP^1 \times\PP^1$ with $2a-2$ vertical lines lying in
$(\Delta_{\infty} \times \PP^1)$, and similarly we let 
$D'_2$ be the union of some other $2b$ horizontal lines in
 $\PP^1 \times\PP^1$ with $2c$ vertical lines lying in
$(\Delta_{\infty} \times \PP^1)$. Let then $D_1, D_2$ be 
respective nearby
smoothings of $D'_1, D'_2$, and let finally let $C$ be the union of
$2$ vertical lines in $(\Delta_0 \times \PP^1)$.

We let $M_1$ be the simple $(\Z/2)^2$ cover of
$ (\Delta_0 \times \PP^1))$ with  branch  curves

$D'_1 \cap (\Delta_0 \times \PP^1))$ and a smoothing of
$ C \cup (D'_1 \cap (\Delta_0 \times \PP^1))$.
 It is rather clear that the roles of $D'_1 , D'_2$ can here be
 freely interchanged by a symmetry in the second $\PP^1$.

Instead , we let $M_2$ be the simple $(\Z/2)^2$ cover of
$ (\Delta_{\infty} \times \PP^1))$ with branch curves
$D_1 \cap (\Delta_{\infty} \times \PP^1))$ and
$D_2 \cap (\Delta_{\infty} \times \PP^1))$.

Another description is as follows: let $X$ be the 
simple $(\Z/2)^3$ cover of $\PP^1 \times\PP^1$
with branch  divisors $D_1, D_2, C$: then $S, S'$ are minimal
 resolutions of the nodal surfaces gotten by dividing $X$ by two
different involutions in the Galois group.

Notice that $\phi^{-1} \circ \phi'$ does not extend to a diffeomorphism
of the simpler manifold $M_1$, and one question is 
whether Floer's theory
could be successfully employed for comparing the two $4$-manifolds.
\end{oss}

\section{Proof of the Theorems }
\label{second}
\begin{proof}{\bf (of Theorem 1.2) }

Let us recall the well known Theorems of Ehresmann and Moser

\begin{teo}(Ehresmann + Moser)
Let $ \pi :  \mathcal {X}\rightarrow T  $ be a proper submersion of
differentiable manifolds with $T$ connected, and assume that we have
a differentiable $2$-form $\omega$ on $\mathcal {X}$ with the property
that

(*) $\forall t \in T$ $ \omega_t : = \omega |_{X_t}$ yields 
a symplectic structure on $X_t$ whose class in $H^2 (X_t, \R)$
 is locally constant on $T$ (e.g., if it lies on $H^2 (X_t, \Z)$).

Then the symplectic manifolds $(X_t, \omega_t)$ are all symplectomorphic.
\end{teo}

Henceforth, applying the lemma to $T: = \Delta - \{0\}$, 
and to the restrictions of the two given families $\mathcal {X}$,
$\mathcal {X}'$, we can for both statements replace $X$ by any
$X_t$ with $t \neq 0$ sufficiently small, and similarly replace $X'$ by
any $X'_{t'}$ with $t' \neq 0$.

In other words, assuming $X_0 = X'_0 \subset \PP^N$, we may assume that
$X$ and $X'$ are both very near to $X_0$.

For each $x_0 \in Sing(X_0)$, $\pi$, resp.$\pi'$, induce germs of 
holomorphic mappings 
$F_{x_0} : \Delta \rightarrow \mathcal{D}_{x_0} : = Def (X_0, x_0)$,
resp. $F'_{x_0}$.

Let $\mathcal{Y}_{x_0} \subset \mathcal{D}_{x_0} \times \PP^N$ be 
the semiuniversal deformation of the germ $(X_0, x_0)$.

For each $ 0 < \epsilon << 1 ,0 < \eta << 1$ 
we consider the family of Milnor 
links 

$$ \mathcal{K}_{\epsilon, \eta} := \mathcal{Y}_{x_0} \cap 
 (B(0, \epsilon ) \times   S (x_0, \eta ) )$$

where $B(0, \epsilon )$ is the ball of radius $\epsilon$ and centre
the point $0 \in \mathcal{D}_{x_0}$ corresponding to $X_0$,
while $S (x_0, \eta ) $ is the sphere in $\PP^N$ with centre $x_0$ and
radius $\eta$ in the Fubini Study metric.

It is well known that, for $\eta << 1$ and $ \epsilon << \eta $,
the family $ \mathcal{K}_{\epsilon, \eta} \rightarrow 
 (B(0, \epsilon ) \cap \mathcal{D}_{x_0}) $ is differentially
trivial (either in the sense of stratified sets, cf. \cite{math},
or, which suffices to us, in the weaker sense that 
when we pull it back through a differentiable
map $\Delta \rightarrow  (B(0, \epsilon ) \cap \mathcal{D}_{x_0}) $
we get a differentiable product).

We use now, to prove statement a), a variant with boundary of
Ehresmann's theorem

\begin{lem}
Let $ \pi :  \mathcal {M} \rightarrow T  $ be a proper submersion of
differentiable manifolds with boundary, such that $T$ is
a ball in $\R^n$, and
assume that we are given a fixed trivialization  $\psi$ of
 a  clased family $  \mathcal {N} \rightarrow T  $ of submanifolds with
boundary. Then we can find a 
trivialization of $ \pi :  \mathcal {M} \rightarrow T  $ which
induces the given trivialization $\psi$.
\end{lem}

\begin{proof}
It suffices to take on $\mathcal {M}$ a Riemannian metric
where the sections $ \psi (p,T) $, for $ p \in  \mathcal {N}$,
are orthogonal to the fibres of $\pi$. Then we use the customary proof
of Ehresmann's theorem, integrating  liftings orthogonal to the fibres of
 standard vector fields on $T$.
\end{proof}

{\bf Proof of a):} we apply lemma 3.2 several times: 

\begin{itemize}
\item
i) first  we apply it in order
to thicken the trivialization of Milnor links to a closed tubular 
neighbourhood
in the semiuniversal deformation, 
\item
ii) then we apply it to the restriction of the families 
$\mathcal {X} \rightarrow \Delta  $, $\mathcal {X}' \rightarrow \Delta  $,
to a ball of radius  $\delta$ where  $\delta$  is so chosen that
$F_{x_0}  (\{ t | \ |t| < \delta \} ) \subset B(0, \epsilon / 2)$
(resp. for $F'_{x_0}$), and to the exterior of the 
balls $B ( x_0, \eta /2)$, so that we  get trivializations
 of the exteriors
to the balls $B ( x_0, \eta /2)$.
\item
iii) we finally use our assumptions that the images of $F'_{x_0}$, resp. 
$F_{x_0}$ land in the same component of $\mathcal{D}_{x_0}$: from 
it follows that there is a holomorphic mapping $ G: \Delta 
\rightarrow \mathcal{D}_{x_0}$ whose image contains the
two points $F_{x_0} (t_0)$, $F'_{x_0}(t'_0)$ and is contained in
$ B(0, \epsilon / 2) \cap \Sigma $ ( $\Sigma$ being as before
 the smoothing locus).

We consider then the pull back to $\Delta$ under $G$ of the family of 
 closed Milnor fibres

$$ \mathcal{M}_{\epsilon, \eta} := \mathcal{Y}_{x_0} \cap 
 (B(0, \epsilon ) \times  \overline{ B (x_0, \eta )} ).$$
To this family we apply again 3.2, in order to obtain a trivialization
of the family of Milnor fibres which extends the given trivialization on
the family of (closed) tubular neighbourhoods of the Milnor links.
\end{itemize}

We are now done, since we obtain the desired diffeomorphism between
$X$ and $X'$ by glueing together (in the intersection with 
$B (x_0, \eta ) - \overline{B (x_0, \eta /2 )}$ ) the two diffeomorphisms
provided by the restrictions of the espective trivializations 
ii) (to the intersection of the complement to 
$\overline{B (x_0, \eta /2 )}$ ) and iii) ( to the intersection
with  $B (x_0, \eta )$): they glue because they both extend
 the trivialization i).

{\bf Proof of b. } The previous construction would allow us to
construct a differentiable family $\mathcal{Z}$ 
over the interval $[-1, 1]$ and
with end fibres $\mathcal{Z}_{-1} \cong X$, resp. 
$\mathcal{Z}_1 \cong X'$. But then it is more cumbersome to show,
by a glueing procedure,
that $\mathcal{Z}$ embeds in $\PP^N \times [-1, 1]$ in such a way that
every fibre inherits a symplectic structure from the Fubini-Study
form ( so that one can then apply Moser's theorem).

We can however more easily construct the desired family $\mathcal{Z}$
of synplectic $4$-manifolds by using the realization of 
symplectic $4$-manifolds  as generic branched covers of $\PP^2$.
This idea goes essentially back to Moisezon (\cite{moi}), we proved the easier
direction (a generic branched cover gives a symplectic $4$-manifold),
the difficult converse result was obtained by Auroux and Katzarkov
in \cite{a-k}: we will refer to this paper, especially to Theorem 3.

By the remark made above, we may concentrate our consideration to the
restriction of the given families $\mathcal{X}$, $\mathcal{X}'$ to
 a disk of radius $\delta << 1$, and assume $X_0 = X'_0$.

It follows that there is a good centre of projection 
 $ L \subset \PP^N$, $L \cong \PP^{N-3}$,
so that the projection $\pi_L : (\PP^N- L) \rightarrow \PP^2 $
is a generic projection for all $X_t, X'_{t'}$.

Let $B_t$ be the branch curve of $\pi_L |_{X_t}$, and similarly
let $B'_{t'}$ be the branch curve of $\pi_L |_{X'_{t'}}$. 
For $t \neq 0 \neq t'$ the corresponding branch curves $B'_{t'}$, 
$B_{t}$ are cuspidal
curves, in the sense that their singularities are only nodes and cusps.
By the cited Theorem 3 of \cite{a-k} it suffices to show that they
are smoothly isotopic in $\PP^2$, or, equivalently, that their associated
Braid Monodromy Factorizations are Hurwitz and conjugation equivalent.

In the terminology of Moisezon, this follows because they provide
the same regeneration of the Braid Monodromy Factorization associated to
$B_0 = B'_0$: we shall try to explain this statement in more detail.

Observe that because the centre of projection $L$ was general,
the curve $B_0$ is a cuspidal curve with the exception of
a finite number of singular points $y_i$ which are the projection 
of exactly one singular point $x_i$ of $X_0$, and of other points where
however the projection $\pi_L | _{X_0}$ is locally invertible.

Consider a small ball $D(y_i,\eta)$ around each such point, and
set $ D := \cup_{i}D(y_i,\eta)$, for $\eta << 1$.

We argue exactly as in part a), that is, we shall prove that 
there is an isotopy between  $B_{t_0}$ and $B'_{t'_0}$ which is
obtained glueing an isotopy in the complement of $D$ and several 
respective isotopies (with fixed boundary) on each $D_i$.

It is clear that for $ |t|, |t'| < \delta << 1$ the 
curves $B_{t} - D$ and 
$B'_{t'} - D$ are isotopic to $B_0 - D$ and the points 
of $B_{t} \cap \partial D$, resp.  $B'_{t'} \cap \partial D$ 
are indeed in a small neighbourhood of $B_0 \cap \partial D$.

On the other hand, consider , for each singular point $x_i$
of $X_0$, the semiuniversal deformation 
$\mathcal{Y}_{x_i} \subset \mathcal{D}_{x_i} \times \PP^N$  
 of the germ $(X_0, x_i)$ and apply the general projection
$\pi_L$ to get a family of deformations 
$\mathcal{B} \rightarrow \mathcal{D}_{x_i}$ of the germ
$(B_0, y_i)$.

By  possibly shrinking  $\mathcal{D}_{x_i}$ ( hence, also $\delta$), 
we get that the family of intersections with the ball 
boundaries, namely, 
 $\mathcal{B}  \cap (\partial D )\times \mathcal{D}_{x_i}$,
is a small deformation of $B_0 \cap \partial D$ and we can find
a trivialization which makes it diffeomorphic to
$ (B_0 \cap \partial D) \times \mathcal{D}_{x_i}$.

We use again the assumption that the two families give two holomorphic
arcs in the  same component of $\mathcal{D}_{x_i}$, whence we use again
the holomorphic mapping $ G: \Delta 
\rightarrow \mathcal{D}_{x_i}$ with image contained in a 
small neighbourhood of the origin and in the smoothing locus,
 and joining the points corresponding to $X_{t_0}$ and  $X'_{t'_0}$.

$G$ induces a family of germs of cuspidal curves which are 
the branch curves of the projections of the Milnor fibres,
and a trivialization of this family gives the required isotopy
of the interior curves $B_{t_0} \cap D_i$ and
$B'_{t'_0} \cap D_i$, which glues together with the given one in the
exterior of $D$ since we may assume as in lemma 3.2 that
 this trivialization extends the one given around the boundary
$\partial D$.

\end{proof}

\begin{proof}{\bf (of Theorem 1.3) }
Let $S$ be the minimal model of a surface of general type.

The assertion is rather clear in the case where the canonical divisor
$ K_S$ is ample.

In fact, let $m$ be such that $m K_S$ is very ample
(any $m\geq 4$ does by Bombieri's theorem, cf. \cite{bom}) thus
the $m$-th pluricanonical map $ \phi_m := \phi_{|mK_S|}$  
is an embedding of $S$ in a projective space $\PP^{P_m-1}$,
where $P_m := dim H^0(\hol_S (m K_S) )$.

We define then  $\omega_m$ as follows:  $\omega_m := \frac{1}{m} \phi_m 
 ^* (FS)$ (where $FS$ 
is the Fubini-Study form $\frac{1}{2 \pi i} \partial 
\overline{\partial} log |z|^2$), whence $\omega_m$ yields 
a symplectic form as desired.

One needs to show that the symplectomorphism class of $(S, \omega_m)$
is independent of $m$. To this purpose, suppose that $n$ is also such
that $ \phi_n$  yields an embedding of $S$: the same holds also for
$nm$, whence it suffices to show that $(S, \omega_m)$ and
$(S, \omega_{mn})$ are symplectomorphic.

To this purpose we use first the well known and easy fact
that the pull back of the Fubini-Study form under the $n$-th
Veronese embedding $v_n$ equals the $n$-th multiple of the Fubini-Study
form. Second,  since $v_n \circ \phi_m$ is a linear projection
of $\phi_{nm}$, by Moser's Theorem follows the desired symplectomorphism.

Assume that $ K_S$ is not ample: then for any $m\geq 5$ (by Bombieri's 
cited theorem) $ \phi_m$  yields an embedding of the canonical model
$X$ of $S$, which is obtained by contracting the finite number of
smooth rational curves with selfintersection number $= -2$ to a finite
number of Rational Double Point singularities. For these, the 
base of the semiuniversal deformation is smooth and yields a 
 smoothing of the singularity.

By Tjurina's theorem (cf. \cite{tju}), $S$ 
is diffeomorphic to any smoothing of $X$:
however we have to be careful because there are many examples 
( cf. e.g. \cite{cat5}) where
$X$ does not admit any global smoothing.

But we observe once more that $S$ is obtained glueing the exterior of
$X - D'$, $D'$ being the union of balls of radius $\eta$ around
the singular points of $X$, with the respective Milnor fibres.

Argueing as in part b) of theorem 1.2 we represent $S$ as generic
branched covering of $\PP^2$ with non holomorphic branch curve,
and we conclude again by theorem $3$ of \cite{a-k} that the canonical 
symplectic structure thus obtained is invariant by smooth deformation
of $S$.
 
\end{proof}

\begin{proof}{\bf  (of Theorem 1.5)}
In \cite{man4} Manetti constructs examples of surfaces $S$, $S'$
of general type which are not deformation equivalent, yet with
the property that
there are flat families of normal surfaces
$  \mathcal {X}\subset \PP^N \times \Delta$ and 
$  \mathcal {X}' \subset \PP^N \times \Delta'$

1) yielding a $\Q$-Gorenstein smoothings of
 the central fibre $X_0 = X_0'$, 

and such that 

2) the  fibres $X_t$, $X_t'$, for $ t, t' \neq 0$ are smooth,
and the canonical divisor of each fibre is ample 

3) there are $t_0$, $t'_0$ with 
$S \cong X_{t_0}$, $S' \cong X'_{t'_0}$.

There exists therefore a positive integer $m$ such that for
each $X_t$ and $X'_{t'}$ the $m$-th multiple of the canonical
(Weil-)divisor is Cartier and very ample, and therefore 
the relative $m$-pluricanonical maps yield two new projective families 
to which 1.2 applies.

By 1.2 and 1.3 follows that $S$ and $S'$, endowed of their canonical 
symplectic structure, are symplectomorphic.

\end{proof}

\section{Application to isotopy of cuspidal plane curves.}

As many authors already observed, results on moduli spaces
of surfaces are strictly intertwined with results on
equisingular families of plane curves. Here is one more specimen

\begin{cor}

 There exist equisingular families of cuspidal algebraic curves
 in $\PP^2$ which are all smoothly isotopic, yet belong to distinct
connected components of the space of cuspidal plane curves of
fixed degree and given number of nodes and cusps. 

\end{cor}
\begin{proof}
As in Theorem 1.5, let us take Manetti's examples of surfaces $S, S'$
and the corresponding cuspidal branch curves $B, B'$ 
for a generic projection
to $\PP^2$. By 1.2, 1.5 the curves $B, B'$ are smoothly isotopic;
however $B, B'$ cannot belong to a connected equisingular family 
of cuspidal plane curves, else the surfaces $S, S'$ would 
be deformation equivalent, a contradiction.
\end{proof}

\bigskip

{\bf Acknowledgements.}
I would like to thank the referee of \cite{cat4} 
for prompting me to write 
down this note, and Anatoly Libgober for suggesting 
to point out the application
to equisingular families of cuspidal plane curves.

\bigskip

\begin{footnotesize}
\noindent
{\bf Note.}
Marco Manetti recently informed me that he was aware of a result
similar to part a) of Theorem 1.2.
\end{footnotesize}

\vfill

\noindent
{\bf Author's address:}

\bigskip

\noindent 
Prof. Fabrizio Catanese\\
Lehrstuhl Mathematik VIII\\
Universit\"at Bayreuth, NWII\\
 D-95440 Bayreuth, Germany

e-mail: Fabrizio.Catanese@uni-bayreuth.de

\end{document}